\documentclass[12pt,reqno]{amsart}

\usepackage{latexsym}
\usepackage{amssymb}
\usepackage[matrix,arrow]{xy}

\author{Florin Ambro} 
\address{
Institute of Mathematics ``Simion Stoilow'' 
of the Romanian Academy, P.O. BOX 1-764, 
RO-014700 Bucharest, Romania.}
\email{florin.ambro@imar.ro}

\newcommand{\Q}{{\mathbb Q}}
\newcommand{\Z}{{\mathbb Z}}

\newcommand{\R}{{\mathbb R}}

\newcommand{\emb}{\operatorname{emb}}

\newcommand{\Int}{\operatorname{int}}

\newcommand{\vol}{\operatorname{vol}}

\theoremstyle{plain}
\newtheorem{thm}{Theorem}[section]

\newtheorem{lem}[thm]{Lemma}

\theoremstyle{definition}

\newtheorem{rem}[thm]{Remark}

\theoremstyle{remark}

\begin{document}

\bibliographystyle{amsalpha+}
\title[Singularities]{On the classification of toric singularities}

\begin{abstract}
For a toric log variety with standard coefficients, 
we show that the minimal log discrepancy at a closed
invariant point bounds the Cartier index of a neighbourhood. 
\end{abstract}

\maketitle


\section{Introduction}
\footnotetext[1]
{
Research partially supported by the grants CEx05-D11-11/04.10.05
and PN-II-ID-PCE-2008-2 cod CNCSIS 2228.
}
\footnotetext[2]{2000 Mathematics Subject Classification. 
Primary: 14B05. Secondary: 14M25.
\clearpage.}

The class of log canonical singularities appears naturally in 
the birational classification of algebraic varieties. 
The main invariants of such singularities are the {\em index} 
and the {\em minimal log discrepancy}, and we expect that these 
two invariants separate the singularities into series. In this 
note we show that the two invariants are equivalent up to finitely 
many values, in the special case of toric singularities.

 To describe the invariants, suppose $P\in X$ is an isolated 
log canonical singularity. Here $X$ is a normal variety. We denote by 
$K$ a canonical Weil divisor of $X$ and suppose $nK\sim 0$ for
a positive integer $n$. Let $n$ be minimal with this property,
called the {\em index} of $P\in X$.
By Hironaka's resolution of singularities, there is a birational 
modification $\mu \colon X'\to X$ such that $X'$ is nonsingular,
$\mu^{-1}(P)=\Sigma$ is a divisor with simple normal crossings,
and $\mu\colon X'\setminus \Sigma\to X\setminus P$ is an isomorphism.
Let $m$ be the smallest multiplicity of the general member of
$|nK_{X'}+n\Sigma|$ along the prime components of $\Sigma$.
Then $\frac{m}{n}$ is independent of the choice of $\mu$
and is called the {\em minimal log discrepancy} of $X$ at $P$, denoted
$a(P;X)$. For the relevance of these invariants to the birational
classification of algebraic varieties, see~\cite{A06}.

Suppose the index is fixed. Then $a(P;X)$ is a rational number with fixed 
denominator. As it is expected that $a(P;X)\le \dim X$, it would
follow that $a(P;X)$ could take only finitely many values. Conversely, 
suppose $a(P;X)$ is fixed. Then, according to Shokurov, we expect that 
the index is bounded. First, this is the analog 
for singularities of the boundedness in terms of volume of canonically 
polarized varieties. Second, there is some evidence for this conjecture. 
A surface germ $P\in X$ with $a(P;X)=0$ has index is $1,2,3,4$ or $6$ 
(Shokurov~\cite{Sh92}). A similar statement holds in dimension $3$ 
(Ishii~\cite{I00}). If $P\in X$ is a terminal $3$-fold singularity, 
then the index is the denominator of $a(P;X)$ (Kawamata~\cite{Sh92}).

It is increasingly becoming clear that in order to classify algebraic 
varieties we must allow not only certain singularities, 
but even certain boundary divisors, to measure ramification. 
These boundaries are crucial in 
the study of singularities and they unify the theories of open and
closed manifolds. In this note we only allow boundaries with so 
called {\em standard coefficients}.

\section{The bound}

We refer to Oda~\cite{O88} for standard notions on toric varieties.
For more details on toric log varieties, see~\cite{A05}.
Let $X$ be an affine toric variety of dimension $d$, 
let $P\in X$ be the unique closed point fixed by the torus.
Let $\{H_\alpha\}$ be the invariant prime divisors of $X$.
Let $B=\sum_\alpha b_\alpha H_\alpha$ be a $\Q$-divisor 
with the following properties:
\begin{itemize}
\item $nK+nB\sim 0$ for some positive integer $n$.
Suppose $n$ is minimal with this property.
\item $\{b_\alpha\}\subset \{\frac{l-1}{l};l\in \Z_{\ge 1}\}\cup\{1\}$.
\end{itemize}

It follows that $(X,B)$ has log canonical singularities and $a(P;X,B)\ge 0$
is a rational number.  

\begin{thm} Let $q$ be the denominator of $a(P;X,B)$. 
Then $n\le c_d q^d$, where $c_d$ is a positive constant depending on 
$d$ only. 

In particular,
if $q$ is fixed then the coefficients of $B$ belong to a finite set.
\end{thm}

\begin{proof} Let $X=T_N\emb(\sigma)$ for a strongly rational polyhedral cone
$\sigma\subset N_\R$. Let $\{e_\alpha\}$ be the primitive points of $N$ on the 
rays of $\sigma$. By assumption, there exists
$\psi\in N^*_\Q$ such that $\langle \psi,e_\alpha\rangle=1-b_\alpha$
for every $\alpha$, and $n\ge 1$ is smallest with $n\psi\in N^*$. 
Since $B$ has one standard coefficient, the sublattice 
$\langle \psi,N\rangle\subset \R$ contains
$1$. Therefore $\langle \psi,N\rangle=\frac{1}{n}\Z$.
The minimal log discrepancy at $P$ is computed as follows:
$$
a(P;X,B)=\min_{e\in N\cap \Int(\sigma)} \langle \psi,e\rangle.
$$
Let $\Lambda=N\cap \psi^\perp$. Choose $e\in N, \langle
\psi,e\rangle=\frac{1}{n}$. Define $\square=\Lambda_\R\cap (\sigma-e)$. 
$e_\alpha=(n-nb_\alpha) (v_\alpha+e)$, where $\{v_\alpha\}$ are the
vertices of $\square$. 
Since $\Lambda\oplus\Z e=N$, we have 
\begin{itemize}
\item $na=\min\{i\ge 1;\Lambda\cap \Int(i\square)\ne \emptyset\}$.
\item $n-nb_\alpha=\min\{i\ge 1; iv_\alpha\in \Lambda\}$ for every vertex $v_\alpha$ of $\square$.
\end{itemize}

Since $B$ is standard, $n\square$ has vertices in $\Lambda$. Denote $na=j$, $S=j\square$. Then 
\begin{itemize}
\item $j=\min\{i\ge 1;\Lambda\cap \Int(\frac{i}{j}S)\ne \emptyset\}$.
\item $S$ has vertices in $\frac{1}{q}\Lambda$.
\end{itemize}

{\em Step 1}: Let $z\in \Lambda\cap \Int(S)$. We may shrink $S$ until
$\{z\}=\frac{1}{q}\Lambda\cap \Int(S)$. Since $S$ has vertices in 
$\frac{1}{q}\Lambda$, it follows by Hensley~\cite{H83} that there
exists a positive constant $\gamma=\gamma_{d-1}$, depending only on $d-1$,
such that $\gamma(S-S)\subset S$.

{\em Step 2}: 
By definition, the cone over $\{j\}\times S$ with vertex $0$ contains no
point of $\Z\times \Lambda$ in its interior. Let $C\subset \R\times
\Lambda_\R$  be the cone over $\{j\}\times (z+\gamma(S-S))$ with vertex $0$. 
Let $C'$ be its reflexion about the lattice point $(j,z)$. Then $P=C\cup C'$ 
is a convex body symmetric about $\{(j,z)\}=\Int(P)\cap \Z\times\Lambda$. 
By Minkowski's first theorem 
(see~\cite{L69} for example), 
$\vol_{\Z\times\Lambda}(P)\le 2^d$. By Lemma~\ref{vo}, this is equivalent to
$$
2\cdot \frac{j}{d}\vol_\Lambda(z+\gamma(S-S))\le 2^d.
$$
By Lemma~\ref{lv}, 
$\vol_\Lambda(z+\gamma(S-S))=\gamma^{d-1}\vol_\Lambda(S-S)\ge
\gamma^{d-1}\frac{2^{d-1}}{(d-1)!q^{d-1}}$. Then 
$
j\le \frac{d!}{\gamma^{d-1}} q^{d-1}.
$
Therefore the claim holds for $c_d=\frac{d!}{\gamma^{d-1}}$.
\end{proof}

\begin{lem}\label{vo}
Let $\square\subset \Lambda_\R, h>0$. Let $C\subset \R\times \Lambda_\R$ be the
cone over $\{h\}\times \square$ with vertex $0$. Then 
$\vol_{\Z\times\Lambda}(C)=\frac{h}{\dim \Lambda+1}\vol_\Lambda(\square)$. 
\end{lem}

\begin{lem}\label{lv}
Let $\square\subset \Lambda^d_\R$ be a lattice convex body. Then 
$\vol_\Lambda(\square-\square)\ge \frac{2^d}{d!}$.
\end{lem}

\begin{proof}
We may assume $\square$ is a simplex with one vertex at the origin, with
vertices $0,v_1,\ldots,v_d$. $\square-\square$ contains the convex hull
$H$ of $\pm v_1,\cdots,\pm v_d$. For $f\in \{0,1\}^d$, denote by $C_f$ the convex
hull of $0,(-1)^{f(1)}v_1,\ldots,(-1)^{f(d)}v_d$. The $C_f$'s cover $H$ and 
have no interior points in common. Each $C_f$ is a lattice convex body, hence 
$\vol_\Lambda(C_f)\ge \frac{1}{d!}$. Their cardinality is $2^d$, hence 
$\vol_\Lambda(\square-\square)\ge \vol_\Lambda(H)=\sum_f \vol_\Lambda(C_f)\ge \frac{2^d}{d!}$.
\end{proof}

\begin{rem}
We may take $c_1=1, c_2=2$. 
\end{rem}


\end{document}